
      \input amssym.def
      \input amssym.tex

      \mag=1200
      \hsize=6.5 true in
      \vsize=8 true in
      \baselineskip=15pt

      \count1=\number\year
      \advance\count1 by -2000
      \def\abc{\number\month/\number\day/0\number\count1}
      \footline={\hss\tenrm -\folio-\ \ \abc \hss}
            


      \def\qed{$\rlap{$\sqcap$}\sqcup$}


    \font\tengothic=eufm10
    \font\sevengothic=eufm7
    \newfam\gothicfam
          \textfont\gothicfam=\tengothic
          \scriptfont\gothicfam=\sevengothic


       \font\tenmsb=msbm10              \font\sevenmsb=msbm7
    \newfam\msbfam
          \textfont\msbfam=\tenmsb
          \scriptfont\msbfam=\sevenmsb
    \def\Bbb#1{{\fam\msbfam #1}}



    \def\move-in{\parshape=1.75true in 5true in}


    \def\proj{{\Bbb P}}


    \def\aff{ {\Bbb A} }






    \hyphenation {Castel-nuovo  Grass-mann-ian}



    \def\twomed{\medskip\medskip}


    \def\ref#1{[{\bf #1}]}

    \vskip 1cm
    \centerline{\bf Secant varieties of Grassmann Varieties}
    \bigskip
    \centerline{\it M.V.Catalisano, A.V.Geramita, A.Gimigliano}

    \bigskip
    \bigskip

    \noindent{\bf Introduction:}
    \medskip

    The problem of determining the dimensions of the higher secant
    varieties of the classically studied projective varieties is a
    problem with a long and interesting history.  The recent work by J.
    Alexander and A. Hirschowitz (see \ref{AH}, for example) completed a
    project that was underway for over 100 years (see \ref{Pa},
    \ref{Te}, and \ref{W}) and confirmed the conjecture that, apart from
    the quadratic Veronese varieties and a (few) well known exceptions,
    all the Veronese varieties have higher secant varieties of the
    expected dimension.

    There has been no comparable success with the case of the Segre
    varieties (nor is there even a compelling conjecture) although there
    is much interest in this question - and not only among geometers.    In 
fact, this particular problem is strongly connected to questions
    in representation theory, coding theory and algebraic complexity
    theory (see our paper \ref{CCG} for some recent results as well as a
    summary of known results, and \ref{BCS} for the connections with
    complexity theory).
    In this paper we investigate this same problem for another family of
    classically studied varieties, namely the Grassmann varieties in
    their Pl\"ucker embeddings.  The dimensions of all the higher secant
    varieties to the Grassmannians of lines in projective space are well
    known (and we give a simple proof of this known result in Section 2)
    but, to the best of our knowledge, little more can be found in the
    classical or modern literature about this problem (although, see the
    comments of Ehrenborg in \ref{E}, and \ref{No} for connections with
    coding theory).

    In this paper we give some general results about the dimensions of
    the higher secant varieties to the Grassmann varieties.  As a
    corollary we can calculate the dimensions of the chordal (or secant
    line) varieties to all the Grassmann varieties.  Our result shows
    that (apart from the known cases) all these secant varieties have
    the expected dimension.

    In Section 3 we discuss our search for deficient secant varieties to
    Grassmannians.  We give new proofs that $G(3,7)$ and $G(3,9)$ have
    deficient secant varieties and also find a new Grassmannian with
    deficient secant varieties, namely $G(4,8)$.
    Finally some words concerning our methods.  Our approach uses, as
    its first step, the fundamental observation of Terracini about
    secant varieties (the so-called {\it Lemma of Terracini}).  This
    Lemma converts the problem of determining the dimension of a
    (higher) secant variety into one of finding the dimension of a
    certain linear space.  The second step in our approach is to
    identify this linear space (via an exterior algebra version of {\it
    apolarity}) with a graded piece of an intersection of homogenous
    \lq\lq fat ideals" in an exterior algebra.  The final step consists
    of calculating the dimension of the appropriate piece of this ideal
    in the exterior algebra.

    \medskip\bigskip \noindent {\bf 1. Preliminaries, Grassmannians and
    Exterior Algebras.}

    \medskip We recall a few elementary facts about exterior algebras
    and Grassmannians.  Let $V$ be a finite dimensional vector space
    over the field $K$ (we will always assume that $char K =0$ and that
    $K$ is algebraically closed) and let $\dim V = n+1$. Then
    $\bigwedge(V)$ denotes the exterior algebra of $V$.  $\bigwedge(V)$
    is a graded (non-commutative) ring such that
    $(\bigwedge(V))_k=\wedge ^kV$.  It is well-known that $\wedge ^k(V)$
    is a finite dimensional vector space and that $\dim \wedge ^k(V)=
    {n+1 \choose k}$.  The elements of $\wedge^k(V)$ are called {\it
    exterior k-vectors} (and sometimes {\it skew-symmetric tensors}). An
    exterior $k$-vector $T$ is said to be {\it decomposable} if we can
    find vectors $v_1, ..., v_k\in V$ such that $T=v_1\wedge ...\wedge
    v_k$. We say that the exterior $k$-vector $T$ has {\it $\wedge$-rank
    $=r$} if it can be written as a sum of $r$ (but no fewer)
    decomposable $k$-exterior vectors.

        It is natural to ask the following two questions:

    \medskip\move-in\noindent $1)$\ What is the least integer $D(k,n+1)$
    so that
    {\it every} exterior vector in $\wedge ^k(V)$ has $\wedge$-rank
    $\leq D(k,n+1)$?

    \medskip\move-in\noindent $2)$\ What is the least integer $E(k,n+1)$
    for which
    there is a dense subset $U \subset \wedge ^k(V)$ (dense in the 
Zariski    topology) such that every exterior vector in $U$ has $\wedge$-rank
    $\leq E(k,n+1)$?
    We will call this number the {\it typical rank} (called the {\it
    essential rank} in [E]) of $\wedge ^k(V)$?

    \medskip Our main focus in this paper will be on Question $2)$.
        \medskip  We first recall the coordinate description of the
    canonical multilinear, alternating map
    $$
    \matrix{ \nu_k : & \underbrace{V\times \cdots \times V}_{k -
    \rm{times}} : &\rightarrow & \wedge ^k (V) \cr \cr
    &(v_1,...,v_k)& \rightarrow & v_1\wedge ...\wedge v_k} \ .
    $$

    Choose an ordered basis $\{ e_0,...,e_n\}$ for $V$.  We identify the
    elements of $\underbrace{V\times \cdots \times V}_{k - \rm{times}}$
    with $k\times (n+1)$ matrices having entries in $K$ by writing
    $v_j=\sum_{i=0}^n a_{i,j}e_i$ and letting $M =(a_{i,j})$ be the
    $k\times (n+1) $ matrix which has the coordinates of the 
vector    $v_j$ as its $j^{\rm{th}}$ row.  We will assume, from now on, that
    we have chosen such an ordered basis for $V$ and that we are making
    the identification above.

    With the basis for $V$ as above, it is standard to choose, as an
    ordered basis for $\wedge^kV$, the set
    $$
    \{ e_{i_1}\wedge \cdots \wedge e_{i_k} \ \mid \ 0 \leq i_1 < \cdots
    < i_k \leq n \}
    $$
    where the subsets $\{ i_1, \ldots , i_k \}$ are ordered
    lexicographically.
    With this choice the map $\nu_k$ can now be written    $$
    \nu_k(v_1, \ldots , v_k) =  \sum _{\{ i_1,...,i_k\} }
    m_{i_1,...,i_k}e_{i_1}\wedge ... \wedge e_{i_k}
    $$
    where $m_{i_1,...,i_k}$ is the $k\times k$-minor of $M$ formed by
    the columns    $i_1, i_2, \ldots , i_k$.

    It is clear from this description that $v_1\wedge ...\wedge v_k \neq
    0$ if and only if the matrix $M$ has rank $k$ i.e. if and only 
if    $\{ v_1, \ldots , v_k \}$ are a linearly independent set of vectors
    if and only if $\langle v_1, \ldots , v_k\rangle$ (the linear space
    spanned by the vectors $v_1, \ldots , v_k$ ) is $k$-dimensional.

    If $\{ v_1^\prime, \ldots , v_k^\prime \}$ is another set of $k$
    vectors from $V$, and $\{v_1^\prime, \cdots , v_k^\prime \}$
    corresponds to the $k\times (n+1)$ matrix  $M^\prime$, then it is a
    standard fact of linear algebra that
    $$
    \langle v_1, \ldots , v_k\rangle = \langle v_1^\prime, \ldots
    ,v_k^\prime\rangle \Leftrightarrow  \hbox{ $\exists$ an invertible
    $k\times k$ matrix $A$ such that } M = AM^\prime .
    $$

    It follows that, in the notation above,
    $$
    m_{i_1, \ldots , i_k} = (\det A)m_{i_1, \ldots , i_k}^\prime \ 
.    \eqno{(\dag)}
    $$

    These last considerations suggest that:

    \medskip\move-in\noindent$1)$\ we consider the map $\nu_k$ as being
    defined on the quotient space of $ \underbrace{V\times \cdots \times
    V}_{k - \rm{times}}$ given by the action of the group $GL_k$ and, in
    particular, only on those elements of the quotient space which have,
    as representative, a set of $k$ linearly independent vectors
    (equivalently, is represented by a $k\times (n+1)$ matrix of rank
    $k$), i.e. we should consider the map $\nu_k$ as being defined on
    the Grassmannian of $k$-dimensional subspaces of $K^{n+1}\simeq V$
    (denoted $G_{k,n+1}$);

    \medskip\move-in\noindent $2)$\ we consider the target space of
    $\nu_k$ as $\proj (\wedge^k(V)) \simeq \proj ^{N_k}$ (where $N_k =
    {n+1\choose k} -1$) and we consider the coordinate functions on
    $\proj ^{N_k}$ to be $z_{i_1, \ldots , i_k}$ where the sets $\{ i_1,
    \ldots , i_k \}$ are chosen as above.

    \medskip\noindent {\bf Remarks:}\
    \noindent $1)$\ None of what we have described above is new, i.e.
    this is the usual description of the Pl\"ucker embedding of the
    Grassmannian of $k$-dimensional subspaces of an $n+1$-dimensional
    vector space into a projective space.  We have insisted on recalling
    all the standard ideas of this embedding because we will be using
    the explicit description of this embedding later and we want to have
    all the notation present for the reader.

    \medskip\noindent $2)$\ Once we choose a basis for $V$ there is a
    natural 1-1 correspondence between the associated bases of
    $\wedge^kV$ and $\wedge^{n+1-k}V$.  It is well known that, if we use
    this corredspondence as an identification, the equations defining
    the Pl\"ucker embeddings of $\wedge^kV$ and $\wedge^{n+1-k}V$ into
    $\proj ^{N_k} = \proj ^{N_{n+1-k}}$ are the same.  Thus, for the
    remainder of this paper we will only consider the Grassmanians
    $G_{k,n+1}$ where $k\leq [(n+1)/2]$.

    \medskip  We will have occasion to use another standard fact from
    linear algebra, namely:  two matrices of size $r\times s$ have the
    same row space if and only if their row reduced echelon forms are equal.

    The way we will use this fact is that we will consider, as
    representatives of the Grassmannian, $G_{k,n+1}$, those matrices of
    size $k\times (n+1)$ which have rank $k$ and which are in
    row-reduced echelon form.  By the comment above, these are precisely
    the matrices which correspond to the different $k$-dimensional
    subspaces of $V$ and on which the map $\nu_k$ is now being defined.

    \medskip\noindent {\bf Definition 1.1:}\ Let $X\subseteq \proj ^N$
    be a closed
    irreducible projective variety; the $s^{th}$ {\it (higher) secant
    variety} of
    $X$ is the closure of the union of all linear spaces spanned by $s$
    points
    of $X$.  We will use $X^s$ to denote the $s^{th}$ secant variety of $X$.
        \medskip There is an \lq\lq expected dimension" for $X^s$: i.e. if
    $\dim X=n$, one \lq\lq expects" that $\dim X^s$ = min$\{N,
    sn+s-1\}$.  This estimate comes from observing that choosing $s$
    points from $X$ gives $sn$ parameters and, having made that choice,
    we are looking at points on a projective space of dimension 
$s-1$.    Thus there are $sn +(s-1)$ parameters, assuming all choices were
    independent.  Of course, the secant varieties to a variety already
    embedded in $\proj ^N$ cannot have a dimension which exceeds 
$N$.    Putting these two things together gives the \lq\lq expected"
    dimension for the variety $X^s$ that we wrote above.
        Since it need not always be the case that the parameters mentioned
    above are independent, it is not always the case that $X^s$ has the
    ``expected dimension".   In such a case we have $\dim X^s <$
    min$\{N, sn+s-1\}$ and we say that $X^s$ is {\it defective}.  A
    measure of this \lq\lq defectiveness" is given by the quantity
    min$\{N, sn+s-1\} - \dim X^s$.
        \bigskip  Let us go back to the embedded Grassmannian,
    $\nu_k(G_{k,n+1}) \subset \proj ^{N_k}$.  Since this variety
    parameterizes the decomposable exterior vectors in $\wedge ^k(V)$,
    its secant varieties $[\nu_k(G_{k,n+1})]^s$ can be viewed as the
    closure of the locus of exterior vectors having $\wedge$-rank $ =
    s$. Hence we have another interpretation of the numbers $E(k,n+1)$
    we introduced above.  More precisely:
        \medskip\noindent {\bf Fact:}\ Let $V$ be a $K$-vector space with
    $\dim_K V= n+1$,  then:
    $$
    E(k,n+1) = \min \{ s \ \vert \  [\nu_k(G_{k,n+1})]^s = \proj ^{N_k} \}.
    $$
    Moreover, information about the dimension of $[\nu_k(G_{k,n+1})]^s
    $, for values of $s$ different from $E(k,n+1)$, will tell us about
    the stratification of $\proj (\wedge ^k V)$ with respect to
    $\wedge$-rank.
        \medskip Before proceeding it is worthwhile to recall that the
    dimensions of all the secant varieties to the Grassmannian of
    2-dimensional subspaces of an $n+1$-dimensional space, i.e. the
    dimensions of $[\nu_2(G_{2,n+1})]^s $, are known (e.g. see \ref{Z}
    ).  Also, some results on $E(k,n+1)$ (from a more algebraic and
    combinatorial point of view) can be found in \ref E.

    \medskip One of the main tools we will use to find the dimensions of
    secant varieties is the famous {\it Lemma di Terracini}.  We recall
    that lemma now.

    \medskip\noindent{\bf Lemma:}\ (Terracini's Lemma) Let $X$ be an
    irreducible variety in $\proj ^n$ and let $P_1, \ldots , P_s$ be
    generic points of $X$.  Let $T_{P_i}$ be the projectivized tangent
    space to $X$ at $P_i$ and denote by $\langle T_{P_1}, \ldots ,
    T_{P_s} \rangle$  the (projective) linear subspace of $\proj ^n$
    spanned by the $T_{P_i}$.  Then,
    $$
    \dim X^s = \dim \langle T_{P_1}, \ldots , T_{P_s} \rangle \ .
    $$

    \medskip In view of Terracini's Lemma, we need a good description of
    the tangent space to a point of $\nu_k(G_{k,n+1}) \subset \proj
    ^{N_k}$.  We do this by studying the differential of the map $\nu_k$
    above.
    Since it is immaterial which point of $\nu_k(G_{k,n+1})$ we
    consider, we will consider the point $\nu_k(M)$, for $M = \pmatrix {
    I_k  & \bf{0} }$, where $I_k$ is the $k\times k$ identity matrix and
    $\bf{0}$ is the $k\times (n+1-k)$ matrix of zeroes.

    The image of this point of the Grassmannian is the point $[1:0:
    \ldots : 0]$ of the Pl\"ucker embedding in $\proj ^{N_k}$.  This is
    a point in the affine piece of $\proj ^{N_k}$ which is the
    complement of the closed set $V(z_{0,1,2,\ldots , k-1})$, i.e. it is
    the point in $\aff ^{N_k}$ whose affine coordinates are $(0,0,\ldots
    ,0)$.

    The points of the Grassmannian, $G_{k,n+1}$, with image in this
    affine piece are all represented by matrices in row reduced echelon
    form of the type $\pmatrix{ I_k & A }$, where $A$ is {\it any}
    matrix of size $k\times (n+1-k)$.  Thus, an affine version of the
    map $\nu_k$ can be described by:
    $$
    \nu_k: \aff ^{k(n+1-k)} \longrightarrow \proj ^{N_k} \backslash
    V(z_{0,1,2,\ldots ,k-1}) \simeq \aff^{N_k} = \aff^{{n+1\choose k} - 
1}    $$
    where, if $A$ is the $k\times (n+1-k)$ matrix which represents a
    point of $\aff^{k(n+1-k)}$, then the coordinates of the image of $A$
    are given by considering all the maximal minors of $\pmatrix{I_k &
    A}$ except the first minor, i.e the minor corresponding to the first
    $k$ columns.

    We now want to compute the linear transformation $d\nu_k(M)$
    explicitly.  Since we know that $\nu_k(G_{k,n+1})$ is a smooth
    variety of dimension $k(n+1-k)$, we are not so interested in the
    rank of this linear transformation (it has to be 
$k(n+1-k)$).    Rather, what we need is an explicit description of the 
vectors in
    the {\it image} of this linear transformation.

    Since $M =\pmatrix{ I_k & {\bf 0 } } $ corresponds to the origin of
    $\aff^{k(n+1-k)}$, the tangent space to $M$ can be thought of as all
    the matrices $B$ of size $k\times (n+1-k)$.  A curve in $\aff
    ^{k(n+1-k)}$ through $M$ with tangent vector $B$ at $M$ is given by
    the line $\lambda B$, where $\lambda \in K$.  The image (under
    $\nu_k$) of a particular point on this line (call it $\lambda_0B$)
    is given by the $k\times k$ minors of the matrix $\pmatrix {I_k &
    \lambda_0B}$ (not the first minor).  The minors which involve all
    but 1 column of $I_k$ give us all the possible $\lambda_0b_{i,j}$,
    where $b_{i,j}$ runs through all the entries of $B$, while those
    minors which involve all but $r$ columns of $I_k$ give us
    $\lambda_0^r(*)$ where $*$ is some $r\times r$ minor of $B$.
    Thus the image of $B$ under the differential $d\nu_k(M)$ is:
    $$
    d\nu_k(M)(B) = lim_{\lambda \rightarrow 0} (1/\lambda)(
    \nu_k(\lambda B) - \nu_k( 0 ) )
    $$

    It is easy to see that, in the limit, we get the entires $b_{i,j}$,
    wherever there was a $\lambda b_{i,j}$, and a 0 wherever there was a
    $\lambda ^r(*)$ for $r > 1$.
    Put another way,  the affine cone over
    $T_{\nu_k(M)}(\nu_k(G_{k,n+1}))$ is:
    $$
     \{ ( \ldots , c_{i_1, \ldots, i_k} , \ldots) \in \aff^{n+1\choose
    k} \mid \matrix { \hbox{where }
    c_{i_1, \ldots , i_k} = 0 \hbox{ if more than one of }\cr \hbox{
    $i_1, \ldots , i_k$ is different from $0,1,\ldots, k-1$ } } \Bigr\} .
    $$
    I.e. it is the vector subspace of $\wedge^k(V)$ generated by all the
    $e_{i_1}\wedge\ldots\wedge e_{i_k}$ where at least $(k-1)$ of the
    $i_j$'s are in $\{0,1,\ldots,k-1\}$.
        \medskip  Recall that there is a natural \lq\lq apolarity" in the
    exterior algebra $\bigwedge (V)$, i.e. there is a perfect pairing
    $$
    \wedge^k(V) \times \wedge ^{n+1-k}(V) \longrightarrow \wedge^{n+1}V
    \simeq K
    $$
    induced by the multiplication in $\bigwedge(V)$.

    Thus, if $Y$ is any subspace of $\wedge^k(V)$, we can associate to
    $Y$, its {\it perpendicular} space $Y^\perp \subseteq
    \wedge^{n+1-k}(V)$ where
    $$
    Y^\perp:= \{ w \in \wedge^{n+1-k}(V) \ \mid \ v\wedge w = 0 \hbox{
    for all } v \in Y \} .
    $$
    Of course, $(Y^\perp )^\perp = Y$ and from standard facts of linear
    algebra we have
    $$
    \dim_K Y^\perp + \dim _K Y = \dim _K \wedge ^k(V) = \dim _K \wedge
    ^{n+1-k}(V) \ .
    $$

    Now, if we let $Y = T_{\nu_k(M)}(\nu_k(G_{k,n+1}))$ (with $M$ as
    above) then
    $$
    Y^\perp = \langle e_{j_1}\wedge \ldots \wedge e_{j_{n+1-k}} \ \mid \
    \hbox{ at least two of } \{ j_1, \ldots , j_{n+1-k} \} \hbox{ are in
    } \{ 0,1,\ldots , k-1 \} \rangle .
    $$

    Put another way,
    $$
    Y^\perp = [(e_0, \ldots , e_{k-1})2]_{n+1-k}  \eqno{(*)}
    $$
    i.e. $Y^\perp$ is the degree $n+1-k$ part of the square of the ideal
    of $\bigwedge(V)$ generated by $e_0, \ldots , e_{k-1}$.

    \medskip Since, quite generally, whenever $W_1, \ldots , W_s$ are
    subspaces of $\wedge^k(V)$ we have that
    $$
    (W_1 + \cdots + W_s)^\perp = W_1^\perp \cap \cdots \cap W_s^\perp
    $$
    we obtain, by applying Terracini's Lemma and the results just
    obtained above, that:

    \medskip\noindent{\bf Proposition 1.2:}\ {\it Let $V$ be a vector
    space of dimension $n+1$ and let}
    $${\cal B}_1 = \{ v_{1,1}, \ldots , v_{k,1} \}, \cdots , {\cal B}_s
    = \{ v_{1,s}, \ldots , v_{k,s} \}
    $$
    {\it be a collection of $s$ sets of $k$ generic vectors in $V$.  Let
    $I_j = (v_{1,j}, \ldots , v_{k,j}) \subset \bigwedge(V)$, $j = 1,
    \ldots ,s$ and let }
    $$
    W = (I_12 \cap \ldots \cap I_s 2)_{n+1-k}
    $$
        {\it Then the dimension of $ [\nu_k(G_{k,n+1})]^s $ is }
    $$ dim_KW^\perp -1 = \biggl[ {n+1\choose k} - dim_KW \biggr] - 1 =
    \biggl[ {n+1\choose k} - 1\biggr] - dim_KW
    $$.  \qed

    \medskip If the subspaces $V_i$, spanned by the ${\cal B}_i$ above,
    are as \lq\lq disjoint as possible", we expect that
    $$
    \dim_KW = max \{ {n+1\choose k} - s[k(n+1-k) + 1], 0 \}
    $$
    i.e. that
    $$
    \dim [\nu_k(G_{k,n+1})]^s = \min\{ {n+1\choose k} -1,
    s[k(n+1-k)]+(s-1) \}
    $$
    which is what we called the {\it expected dimension } of
    $[\nu_k(G_{k,n+1})]^s $.
    I.e. if $W$ has the expected dimension then so does
    $[\nu_k(G_{k,n+1})]^s$.  Moreover, if $W$ has more than the {\it
    expected dimension} then the difference is precisely what we called
    the {\it defectiveness} of $[\nu_k(G_{k,n+1})]^s$.

    \medskip\noindent{\bf Remark:}\ Notice that as long as $ks \leq n+1$
    we can choose the vectors $v_{i,j}$ to be part of the basis $\{ e_0,
    \ldots , e_n \}$ of $V$.  It is precisely this case we will consider
    in the next section.

    \bigskip\noindent {\bf 2. The ``monomial" case.}
        \medskip In this section we will suppose that $ks\leq n+1$ and let
    $$
    W=[(e_0,...,e_{k-1})2\cap ...\cap (e_{ks-k},...,e_{ks-1}) 2]_{n+1-k}.
    $$
    By analogy with the case of ideals in the symmetric algebra of a
    free module, we call call such ideals "monomial ideals" of the
    exterior algebra.
    Thus, we can view $W$ as the degree $n+1-k$ part of a monomial ideal
    (the intersection of $s$ squares of monomial ideals).

    \medskip We have the following theorem:
        \medskip\noindent {\bf Theorem 2.1:}\ {\it Let V,n,k be as in the
    previous section. Then:}
        {\it i)  if $k=2$ then  $[\nu_2(G_{2,n+1})]^s$ is defective for } $s
    < E(2,n+1)=\lfloor {n+1 \over 2}\rfloor$ {\it with defectiveness }
    $2s(s-1)$;
        {\it ii)  while if $k\geq 3$ and $ks\leq n+1$, then
    $[\nu_k(G_{k,n+1})]^s$ has    the expected dimension.}
        \medskip

    \noindent {\it Proof:} The case $k=2$ is known (e.g. see \ref E or
    also \ref Z)
    but we give the proof here for the sake of completeness and also
    because the    ``monomial ideal" approach makes it quite easy.
        First, assume that $2s\leq n+1$. We have to consider the vector space
    $$
    W=[(e_0,e_{1})2\cap (e_2,e_3) 2\cap ...\cap
    (e_{2s-2},e_{2s-1})2]_{n-1}
    $$

    Once we note that $(e_i,e_j)2 = (e_i\wedge e_j)$ we have that
    $$
    W = (e_0\wedge e_1 \wedge \cdots \wedge e_{2s-2} \wedge e_{2s-1})_{n-1}
    $$
    Thus, $W$ will trivially be $=\{ 0\}$ if and only if $n-1 < 
2s$.    This immediately gives that $$E(2,n+1)=\lfloor {n+1 \over 2}\rfloor
    . $$
    When $n-1 \geq 2s$, we get that a basis for $W$ is given by the
    decomposable exterior vectors of the form:
    $$
    e_0\wedge e_1\wedge ... \wedge e_{2s-1}\wedge e_{\alpha _1}\wedge ...
    \wedge e_{\alpha _t},
    $$     where $t=n-1-2s$ and $\alpha_1, \ldots , \alpha_t$ can be any $t$
    elements in $\{ 2s, \ldots , n \}$.  So,    $$
    \dim_K W = \dim _K ( \wedge ^t \langle e_{2s},...,e_{n}\rangle ) =
    {n-2s+1 \choose n-2s-1} = {n-2s+1 \choose 2}.
    $$
    A simple computation shows that this is not the expected dimension for
    $W$ (the expected dimension is ${n+1\choose 2}-s(2n-1)$ in this
    case) and that    the defectiveness is $\delta = 2s 
2-2s=2s(s-1)$.  This completes the
    proof of $i)$.

    \medskip  As for $ii)$, let $k\geq 3$ and $ks\leq n+1$. Recall that
    we may always suppose that $k \leq (n+1)/2$ (see Remark 2 before
    Definition 1.1).

    \medskip\noindent{\bf Case 1:}\ Suppose that $2s > n+1 - k$.

    In order for there to be a monomial $e_{i_1}\wedge \ldots \wedge
    e_{i_{n+1-k}}$ in $W$ we must have at least two of $\{i_1, \ldots ,
    i_{n+1-k} \}$ from each of the $s$ subsets
    $$
    \{ 0, \ldots , k-1\}, \ \{k, \ldots , 2k-1 \}, \ \cdots \ ,\{ks-k,
    \ldots , ks-1 \} .
    $$
    So, if $2s > n+1-k$ this will automatically give that $W = 0$.  It
    remains to check that, under the given hypothesis on $s,n$ and $k$,
    the {\it expected dimension} of $W$ is also 0.

    Observe that in order to have: $k\geq 3, \ ks \leq n+1$ and $2s >
    n+1-k$ we must have:
    $$
    { {n+1-k}\over 2} < s \leq {{n+1}\over k} . \eqno{(\dag \dag)}
    $$
    This implies that
    $$
    k2 - k(n+1) + 2(n+1) > 0 \ \eqno{(\ddag)}
    $$

    The discriminant of the quadratic expression in $(\ddag)$ is $\Delta
    = (n+1)(n-7)$, so for $n<7$, the quadratic expression has the same
    sign for every $k$.  Since the coefficient of $k 2$ is positive this
    means that $(\ddag)$ is always satisfied for such $n$.

    Now, the conditions: $k\geq 3$, $n\leq 6$, $ks \leq 7$,  and $2k
    \leq 7$ has only one solution:  $k=3,\ s=2, \ n=5$.  A quick check
    shows that for these parameters the {\it expected dimension} for $W$
    is indeed $0$.  This gives $E(3,6) = 2$ and $\nu_3(G_{3,6}) 2 =
    \proj ^{14}$ as required.

    When $n=7$, the only possibility for $k$ is 3 and there is no $s$
    satisfying $(\dag\dag)$.

    So, suppose that $n \geq 8$.  The quadratic equation associated to
    $(\ddag)$ then has two distinct roots
    $$
    r_1 = {{n+1}\over 2} - { {\sqrt{\Delta}}\over 2} < {{n+1}\over 2} +
    { {\sqrt{\Delta}}\over 2} = r_2
    $$
    Thus $(\ddag)$ can be satisfied only for $k\leq r_1$ and $k \geq
    r_2$.  But, $r_2 > {{n+1} \over 2}$, so we need only consider $k
    \leq r_1$.

    Notice that, for $n>8$, $  {{n+1}\over 2} - { {\sqrt{\Delta}}\over
    2} < 3$, so only case left to consider is $n=8$ and $k=3$ and there
    is no $s$ satisfying $(\dag\dag)$ for those values of $n$ and $k$.

    This completes Case 1.

    \medskip\noindent{\bf Case 2:}\ $2s \leq n+1-k$.

    In this case $W$ is certainly $\neq 0$.  Since $W$ is the degree
    $n+1-k$ part of a monomial ideal, to compute $\dim_KW$ it is enough
    to count all the monomials of degree $n+1-k$ which are NOT in $W$.

    These are the monomials $e_{\alpha_1}\wedge \ldots \wedge
    e_{\alpha_{n+1-k}}$ with the property that, when we choose
    $\alpha_1, \ldots , \alpha_{n+1-k}$ from $\{ 0 , 1, \ldots , n \}$
    we must make sure that from {\bf at least} one of the $s$ subsets
    $$
    \{0, \ldots , k-1 \}, \{k, \ldots , 2k-1 \}, \cdots , \{(s-1)k,
    \ldots , sk-1 \}
    $$
    we have chosen either nothing or one index.

    If we concentrate at first (say) on the subset $\{ 0, \ldots , k-1
    \}$, then we need to choose $n+1-k$ elements from $\{0, \ldots , n
    \}$ such that all $n+1-k$ are outside $\{ 0, \ldots , k-1 \}$ or at
    most one is in $\{0, , \ldots , k-1 \}$

    Since there are exactly $n+1-k$ elements outside $\{0, \ldots , k-1
    \}$, we have only 1 choice if we choose nothing from $\{0, \ldots ,
    k-1 \}$.  However, there are $n+1-k$ subsets of $\{k, \ldots , n \}$
    consisting of $(n+1-k)-1$ elements and to those we can add any one
    of the $k$ elements in $\{0, \ldots , k-1 \}$.  This gives us a
    total of $1 + k(n+1-k)$ choices.  Since we can do this for each of
    the $s$ subsets of $k$ elements above, we get a total of $s[k(n+1-k)
    +1]$ choices.  I.e. there are exactly $s[k(n+1-k) +1]$ monomial of
    degree $n+1-k$ outside $W$.
    Since this is precisely the {\it expected codimension} of $W$, the
    theorem is proved.

    \medskip The previous result has the following immediate corollary:
        \twomed\noindent{\bf Corollary 2.2:}\ {\it  Let $V,\ n, \ k$ be as
    above but suppose that   $k\geq 3$. Then all the chordal varieties,
    $ \nu_k(G_{k,n+1})2$ have the expected dimension.}
             \medskip\medskip\noindent  {\bf 3. Some Final Remarks }

    Since there appears to be so little in the literature about secant
    varieties to Grassmannians (but lots of folklore) we would like to
    collect some scattered results we have found and record them in this
    section.  This will include an (apparently) new example of a
    deficient variety.

    First, let $m,\ n,\ d$ be three positive integers such that $m = n +
    d - 1$.  Let $V$, $W$ be vector spaces over $K$ such that $\dim _KV
    = n$ and $\dim _KW = m$.  Let $\phi_d: \proj ^{n-1} \longrightarrow
    \proj ^t$ be the Veronese embedding (and so $t = {m\choose d} -
    1$).  The Pl\"ucker embedding of $\wedge ^dW$ is also in $\proj ^t$
    since, as Ehrenborg pointed out in \ref{E}, there is a 1-1
    correspondence beetween $d$-subsets of an $m$-set and
    $d$-multisubsets of an $n$-set.

    This suggests that in looking for deficient Grassmannians we
    consider all the deficient Veronese varieties.  We do that now.

    \medskip\noindent{\bf Case 1:} \ All the quadratic Veronese
    varieties are defective, i.e. when $d=2$ we can choose $n$
    arbitrarily.  What about the corresponding Grassmannians?  For $m =
    4,\ 5$ the Grassmannian of lines is not defective, for trivial
    reasons.  However, it is the case that all the other Grassmannians
    of lines are defective.  (See Theorem 2.1 i) ).  This is hopeful.

    \medskip\noindent{\bf Case 2}:\ There are exactly four other
    defective Veronese varieties.  They correspond to:

    \move-in\noindent $i)$\ $d=4$, $n=3$ and hence $m=6$;

    \move-in\noindent $ii)$\ $d=4$, $n=4$ and hence $m = 7$;

    \move-in\noindent $iii)$\ $d=3$, $n=5$ and hence $m=7$;

    \move-in\noindent $iv)$\ $d=4$, $n=5$ and hence $m=8$.

    Now, $i)$ suggests we look at $\nu_4(G_{4,6})$.  Since $\wedge^4K 6
    \simeq \wedge ^2K6$, this is indeed defective, it is (again) a
    Grassmannian of lines.

     >From $ii)$  we consider $\nu_4(G_{4,7})$. Since $\wedge 4 K7
    \simeq \wedge3 K 7$ we should check if $\nu_3(G_{3,7})$ is
    defective.  The secant plane variety should fill $\proj ^{34}$ but
    it does not, as we will now show.

    By our method, we have to determine $\dim Z$, where
    $$
    Z= [(e_0,e_1,e_2)2\cap (e_3,e_4,e_5) 2 \cap (e_6,v,w)2]_4
    $$
    where $\{ e_0,...,e_6\}$ is a basis of $W$ and $v$, $w$ are generic
    vectors
    in $W$.
    \par
    We can actually suppose that $v$,$w$ $\in \langle e_0,...,e_5
    \rangle$, in fact if
    $v=e_6+v_1$ and $w=e_6+w_1$ we get:
    $$
    (e_6,v,w)2 = (e_6\wedge v, e_6\wedge w, v\wedge w) =
    (e_6\wedge v_1, e_6\wedge w_1, v_1\wedge w_1) = (e_6,v_1,w_1) 2.
    $$
    Notice that:
    $$
    Z' = [(e_0,e_1,e_2)2\cap (e_3,e_4,e_5) 2]_4 =
    \langle e_{i_1}\wedge e_{i_2} \wedge e_{i_3} \wedge e_{i_4}\rangle 
_4,    $$
    $$
    {\rm with} \qquad  i_1<i_2\in \{0,1,2\},
    \quad i_3<i_4 \in \{3,4,5\} .
    $$
    Thus, if there is something in $Z$ it must be of the form $v\wedge
    w\wedge \Gamma$,
    $\Gamma \in \wedge ^2\langle e_0,...,e_5\rangle$. Actually, by the
    genericity of $v,w$ we
    can suppose that $\langle e_0,...,e_5\rangle = \langle
    v,w,e_1,...,e_4\rangle$, hence we can consider
    $\Gamma \in \langle e_1,...,e_4\rangle$.
        We can even forther suppose that $v=e_0+...+e_5$ and $w =
    \sum_{i=0}^5 b_ie_i$ and then it is enough to insure that for each
    of the 6 monomials, $m$, in $\Gamma$, the summands of $v\wedge w
    \wedge m$ which are not in $Z^\prime$ are all $= 0$.
        This gives us 6 linear equations in the $b$'s whose coefficient
    matrix is:
        $$
    \pmatrix {b_3-b_0&b_0-b_2&0&b_1-b_0&0&0\cr
    b_4-b_0&0&b_0-b_2&0&b_1-b_0&0\cr
    b_5-b_0&0&0&0&0&0\cr
    0&0&0&0&0&b_5-b_0\cr
    0&b_5-b_4&b_3-b_5&0&0&b_5-b_1\cr
    0&0&0&b_5-b_4&b_3-b_5&b_5-b_2}
    $$
    This is a matrix whose rank is $=5$ and hence we get  $\dim Z=1$,
    i.e. $\nu_3(G_{3,7})$ is defective for secant $\proj ^2$'s, with
    defectivity 1.

    \medskip\move-in\noindent {\bf Note:}\ Apparently this example is
    well-known, as we recently learned from J. Landsberg (private
    communication).  We were unable to find a reference to it in the
    literature.

    \medskip From $iii)$ we are again brought to consideer
    $\nu_3(G_{3,7})$, which we have just done.

    \medskip From,$iv)$ we should consider $\nu_4(G_{4,8})$.  We now
    show that this is indeed a defective variety.

    The Grassmannian $\nu_4(G_{4,8})$  lies in
    $\proj^{69}$.  Since $\dim\nu_4( G_{4,8}) =16$ we have, by our
    Corollary 2.2, that
    $\dim(\nu_4( G_{4,8}))2=33$, the expected dimension.  One expects that
    $\dim (\nu_4(G_{4,8}))3= 50$, but we will show that $
    (\nu_4(G_{4,8}))3=49$ instead.

    By what we have seen before, it will be enough to prove that $\dim_K
    W=20$, where if $H_1, H_2$ and $H_3$ are three generic subspaces of
    $K8$, each of dimension 4, with bases $\{v_{i1}, v_{i2}, v_{i3},
    v_{i4} \}$, $i = 1,2,3$, then
    $$
    W = [(v_{11}, v_{12}, v_{13}, v_{14})2 \cap (v_{21}, v_{22},
    v_{23}, v_{24})2 \cap (v_{31}, v_{32}, v_{33}, v_{34})2]_4
    $$
    in the exterior algebra $\wedge(K8)$.

    We can assume that $H_1 = \langle e_0, e_1, e_2,e_3 \rangle$, where
    the $e_i$ are part of a standard basis for $K 8$.  Consider $H_3
    \cap \langle H_2, e_i \rangle$.  This is a one dimensional subspace
    of $H_3$ which we'll denote by $\langle w_i \rangle$  By the
    genericity of the subspaces we may suppose that $\{ w_0, w_1, w_2,
    w_3 \}$ are a basis for $H_3$.  But, $w_i = e_i + u_i$ for $u_i \in
    H_2$ and again, using the genericity, we can assume that $\{ u_0,
    u_1, u_2, u_3 \}$ are a basis for $H_2$.

    Putting all this together we can assume, without loss of generality,
    that
    $$
    H_1 = \langle e_0, e_1, e_2, e_3 \rangle , \ H_2 = \langle e_4, e_5,
    e_6, e_7 \rangle \hbox{ with } \{ e_0, \ldots , e_7 \} \hbox { a
    basis for } K8
    $$
    and
    $$
    H_3 = \langle e_0 + e_4, e_1 + e_5, e_2 + e_6, e_3 + e_7 \rangle .
    $$

    Now a simple calculation using {\it Macaulay 2} shows that $\dim_KW$
    is indeed 20 and we are done.  In fact, a simple hand check shows
    that the 20 forms in this space are:
    $$
    e_0\wedge e_1\wedge e_4\wedge e_5; \
    e_0\wedge e_2\wedge e_4\wedge e_6; \
    e_0\wedge e_3\wedge e_4\wedge e_7;
    $$
    $$    e_1\wedge e_2\wedge e_5\wedge e_6; \
    e_1\wedge e_3\wedge e_5\wedge e_7; \
    e_2\wedge e_3\wedge e_6\wedge e_7,
    $$
    the 12 forms
    $$
    e_i\wedge e_{i+4}\wedge(e_j\wedge e_{k+4}\pm e_k\wedge e_{j+4}); \quad
    i\neq j \neq k, \quad i,j,k\in \{0,1,2,3\}
    $$
    and the two forms
    $$
    (e_0\wedge e_5\pm e_1\wedge e_4)\wedge (e_2\wedge e_7\pm e_3\wedge
    e_6);
    \qquad    (e_0\wedge e_6\pm e_2\wedge e_4)\wedge (e_1\wedge e_6\pm 
e_3\wedge
    e_5).
    $$

    This shows that$(\nu_4(G_{4,8}))3$  is 1-defective.
    One wonders about the dimensions of the other secant varieties of
    $\nu_4(G_{4,8})$.  A calculation (using  {\it Macaulay 2}) shows
    that $(\nu_4(G_{4,8}))5 = \proj ^{69}$.
    It is not hard to show that if we fix three general points $P$, $Q$,
    and $R$ on the Grassmannian in $\proj ^{69}$ then the linear system
    of hyperplanes which contain the tangent spaces $T_P$ and $T_Q$ and
    which also contain $R$ have exactly a fixed line in the tangent
    plane $T_R$.  (This is the geometric reason why the secant planes to
    $\nu_4(G_{4,8})$ are 1-deficient.)

    Hence, when we take a fourth point $S$ on $\nu_4(G_{4,8})$ and
    consider the linear system of hyperplanes on $\proj ^{69}$ which
    contains $T_P$, $T_Q$, and $T_R$ and which also contain S, that
    system contains three fixed lines in $T_S$.  We expect those lines
    to be linearly independent and hence that the 4-secants to
    $\nu_4(G_{4,8})$ are 4-defective.  Calculations with {\it Macaulay
    2} seem to confirm that expectation.

    \twomed It would be tempting, at this point, to conjecture that
    inasmuch as we have exhausted the list of defective Veronese
    varieties then we have also exhausted the list of defective
    Grassmannians!  Indeed, we hoped that this might be so, but J.
    Landsberg informed us that he had a communication from M.
    Catalano-Johnson asserting that $\nu_3(G_{3,9})$ is also defective.

    Since we couldn't find that example in the literature we provide a
    proof now.
    Notice that, in view of Theorem 2.1 ii, the space of secant $\proj
    ^2$'s does have the correct dimension.  So, we will now show that
    $(\nu_3(G_{3,9}))4$ has dimension 73 instead of the expected
    dimension 75.

    The argument follows the same lines we used to find the defectivity
    of $\nu_4(G_{4,8})3$.  Following the discussion in $\S$ 2 we need
    to find $\dim_KW$, where if $H_i, \ i = 1, \ldots , 4$ are 4
    generics 3-dimensional subspaces of $K9 $ and a basis for $H_i$ is
    $\{v_{i1}, v_{i2}, v_{i3}\} $, then
    $$
    W = \bigl[ \bigcap_{i=1}^4 (v_{i1}, v_{i2}, v_{i3}) 2 \bigr] _6 .
    $$

    It is easy to see that, with no loss of generality, we can assume
    the four subspaces are:
    $$
    H_1 = \langle e_1, e_2, e_3 \rangle, \ H_2 = \langle e_4, e_5, e_6
    \rangle,\ H_3 = \langle e_7, e_8, e_9 \rangle
    $$
    and
    $$
    H_4 = \langle e_1 + e_4 + e_7, e_2 + e_5 + e_8, e_3 + e_6 + e_9 \rangle.
    $$
    Using the exterior algebra routines in {\it Macaulay 2} we find that
    $\dim_KW = 10$ and so the dimension of $\nu_4(G_(4,8)) 3$ is 73, as
    stated.

    These last two examples suggested that we should check
    $\nu_3(G_{3,12})$ and $\nu_4(G_{4,12})$ as well.  We have verified
    that $\nu_3(G_{3,12})5$ and $\nu_4(G_{4,12}) 4$ are {\bf not}
    defective.

    \vfill\eject

    \centerline {{\bf REFERENCES}}
    \par
    \medskip
    \medskip

    \medskip
    \noindent [{\bf AH}]: J.Alexander, A.Hirschovitz.
    {\it Polynomial interpolation in several variables.}
    J. of Alg. Geom. {\bf 4} (1995). 201-222.

    \medskip
    \noindent [{\bf BCS}]: P. B\"{u}rgisser, M. Clausen, M.A.
    Shokrollahi, {\it Algebraic Complexity Theory}, Vol. 315, Grund.
    der Math. Wiss., Springer, 1997

    \medskip
    \noindent [{\bf CGG}]: M.V.Catalisano, A.V.Geramita, A.Gimigliano. {\it
    Rank of tensors, Secant Varieties of Segre Varieties and Fat Points.}
     Preprint, To appear inn Lin. Algebra and Appl.

    \medskip
    \noindent [{\bf E}]: R.Ehrenborg. {\it On Apolarity and Generic
    Canonical Forms.} J. of Algebra {\bf 213} (1999). 167-194.

    \medskip\noindent [{\bf No}]: D. Yu. Nogin. {\it Spectrum of codes
    associated with the Grassmannian G(3,9).} \ Problems of Information
    Transmission, {\bf 33} (1997). 114-123.

    \medskip\noindent [{\bf Pa}]: F.Palatini. {\it Sulle variet\`a
    algebriche per le
    quali sono di dimensione minore} {\it dell' ordinario, senza
    riempire lo
    spazio ambiente, una o alcuna delle variet\`a} {\it formate da spazi
    seganti.}
    Atti Accad. Torino Cl. Scienze Mat. Fis. Nat. {\bf 44} (1909). 362-375.


    \medskip\noindent [{\bf Te}]: A.Terracini. {\it Sulle} $V_k$ {\it
    per cui la
    variet\`a degli} $S_h$ $(h+1)${\it -seganti ha dimensione minore
    dell'ordinario.} Rend. Circ. Mat. Palermo {\bf 31} (1911). 392-396.

    \medskip\noindent [{\bf W}]: W.Wakeford. {\it On canonical forms.}
    Proc. London Math. Soc. {\bf 18} (1919/20). 403-410.

    \medskip\noindent [{\bf Z}]: F.L.Zak. {\it Tangents and Secants of
    Algebraic
    Varieties.} Translations of Math. Monographs, vol. {\bf 127} AMS.
    Providence
    (1993).

    \medskip
    \bigskip
    {\it M.V.Catalisano, Dip. Matematica, Univ. di Genova, Italy.}

    {\it e-mail:\ catalisa@dima.unige.it }
    \medskip

    {\it A.V.Geramita, Dept. Math. and Stats. Queens' Univ. Kingston,
    Ont., Canada}
    {\it and Dip. di Matematica, Univ. di Genova. Italy.}

    {\it e-mail:\  geramita@dima.unige.it ; tony@mast.queensu.ca }

    \medskip
    {\it A.Gimigliano, Dip. di Matematica and CIRAM, Univ. di Bologna,
    Italy.}

    {\it e-mail:\ gimiglia@dm.unibo.it }

    \end